\documentstyle[12pt]{article}

\begin{document}

\baselineskip=18pt

\title{{ Properties at infinity of diffusion semigroups \\
	and stochastic flows via  weak uniform covers}}
\author{  Xue-Mei LI}
\date{ }
\maketitle

\newcommand{\A}{{\bf \cal A}}
\newcommand{\B}{{ \bf \cal B }}
\newcommand{\C}{{\cal C}}
\newcommand{\F}{{\cal F}}
\newcommand{\G}{{\cal G}}
\newcommand{\h}{{\cal H}}
\newcommand{\K}{{\cal K}}
\newcommand{\half}{{  {1\over 2}  }}
\newcommand{\heatsemif}{{ {\rm e}^{ \half t\triangle^{h,1}}   }}
\newcommand{\heatsemi}{{ {\rm e}^{\half t \triangle^{h}}   }}

\newtheorem{theorem}{Theorem}[section]
\newtheorem{proposition}[theorem]{Proposition}
\newtheorem{lemma}[theorem]{Lemma}
\newtheorem{corollary}[theorem]{Corollary}
\newtheorem{definition}{Definition}[section]

\def\limsup{\mathop{\overline{\rm lim}}}
\def\liminf{\mathop{\underline{\rm lim}}}
\def\exp{{\rm e}}

\begin{abstract}
A unified treatment is given of some results of H. Donnelly-P. Li and 
L. Schwartz concerning the behaviour of heat semigroups  on open manifolds 
with given compactifications, on one hand, and the relationship with the
 behaviour at infinity of solutions of related stochastic differential
 equations on the other. A principal tool is the use of certain covers of 
the manifold: which also gives a non-explosion test. As a corollary we
 obtain known results about the behaviour of   Brownian motions on a
 complete Riemannian manifold   with Ricci curvature decaying at most
 quadratically in the distance function.

\end{abstract}


\section{Introduction}

Let $M$ be a manifold, $\{\Omega,{\cal F},{\cal F}_t, P\}$ a filtered
 probability space.  A diffusion process on $M$  is a stochastic process
 which is path continuous  and strongly Markov. These diffusions are often
 given by solutions to  stochastic  differential equations of the following
 type:  
	\begin{equation}\label{eq: basic}
dx_t=X(x_t)\circ dB_t+A(x_t)dt 
\end{equation}
Here $B_t$ is a $R^m$ valued Brownian motion,   $X$ is $C^2$ from
 $ R^m\times M$ to the tangent bundle $TM$  with $X(x)$: $R^m$ $\to T_xM$ a
 linear map for each $x$ in  $M$,  $A$ is a $C^1$ vector field on $M$. Let
  $u$  be an $M$-valued random variable independent of $\F_0$.  Denote by 
 $F_t(u)$ the solution starting from $u$, $\xi(u)$ the explosion time. Then
 $\{F_t(x)\}$ is a diffusion for each $x\in M$. Associated with $F_t$ 
there are  also the probabilistic semigroup $P_t$ and the corresponding 
infinitesimal generator $\A$. 

\bigskip

\noindent
{\bf I.  Main results: } 
 The main aim \cite{LI89}  of this article is to give 
unified treatment to some of the results from H. Donnelly-P. Li and L.
 Schwartz. It gives the first a probabilistic interpretation and extends part
 of the latter. We first introduce  weak uniform covers in an analogous way
 to  uniform covers, which gives a nonexplosion test by using estimates on
 exit times of the diffusion considered. As a corollary this gives the known
 result on nonexplosion of a Brownian motion on a complete Riemannian manifold
 with Ricci curvature decaying at most quadratically in the distance function
\cite{IC82II}.

One interesting example is that a solution to a stochastic differential equation on $R^n$ whose coefficients have linear growth has no explosion and has the $C_0$ property. Notice under this condition,
its associated generator has quadratic growth. On the other hand let $M=R^n$,
 and let $L$ be an elliptic  differential operator :
$$L=\sum_{i,j} a_{ij}{{\partial^2}\over {\partial x_i\partial x_j}}
       +\sum_i b_i{\partial \over {\partial x_i}}$$
where $a_{ij}$ and $b_i$ are $C^2$. Let $(s_{ij})$ be the positive square root of 
the matrix $(a_{ij})$. Let
  $X^i=\sum_j s_{ij}{\partial \over \partial x_j}$,
$A=\sum_j b_j {\partial \over {\partial x_j}}$.  Then 
the s.d.e. defined by:
$$(It\widehat o)  \hskip 20pt    dx_t=\sum_iX^i(x_t)dB_t^i +A(x_t)dt$$
has generator $L$. Furthermore if $|(a_{ij})|$ has quadratic growth and $b_i$ has 
linear growth,  then both $X$ and $A$ in the s.d.e. above have linear growth.
In this case any solution $u_t$ to the following partial differential equation:

$${{\partial u_t}\over {\partial t}} =Lu_t $$ 

\noindent 
satisfies : $u_t\in C_0(M)$, if $u_0\in C_0(M)$ (see next part). 

\bigskip

\noindent
{\bf II.  Preliminaries}

A diffusion process is said to be a $C_0$ diffusion if its semigroup leaves
 invariant $C_0(M)$, the space of continuous functions vanishing at infinity,
 in which case the   semigroup is said to have the $C_0$ property. A
 Riemannian manifold is said to be  stochastically complete
if the Brownian motion on it is complete, it is also said to have the $C_0$
 property if the Brownian motion on it does. A Brownian motion is, by
 definition, a path   continuous strong Markov process with generator
 ${1\over 2 } \triangle$, where  $\triangle$ denotes the Laplacian Beltrami
 operator.  Equivalently a manifold is stochastically complete if the minimal
 heat semigroup  satisfies: $P_t1\equiv 1$.

One example of a  Riemannian manifold which is stochastically complete is a
 complete manifold  with finite volume. See Gaffney \cite {GA59}.  More
 generally a complete  Riemannian manifold with Ricci curvature bounded from
 below  is stochastically   complete and 
has the $C_0$ property as proved by Yau \cite {Yau78}. See also 
 Ichihara \cite{IC82II}, Dodziuk \cite{Dodz}, Karp-P. Li \cite{KA-LI},
 Bakry \cite{BA86},  Grigor\'yan \cite{GR87}, Hsu\cite{Hsu89},
 Davies \cite{DA89}, and Takeda \cite{TA91}  for further discussions in terms
 of volume growth and bounds on Ricci curvature.  For discussions on the
 behaviour at infinity of diffusion processes, and the $C_0$ property, we
 refer the reader to Azencott \cite{AZ74} and Elworthy \cite{EL82}. 

Those papers above are on a Riemannian manifold except for the last reference, 
for a manifold without a Riemannian structure, Elworthy \cite{ELbook}
 following It\^ o \cite {ITO50}  showed that the diffusion  solution
 to (1) does not explode if there is a uniform cover for the coefficients of
 the equation. See also Clarke\cite{Clark73}. In particular this shows that
 the s.d.e (1) does not explode on a
 compact manifold if the coefficients are reasonably smooth. See  
\cite{CA-EL83}, \cite{ELflour}.  To apply this method to check whether a 
Riemannian manifold is stochastically complete, we usually construct a 
stochastic differential equation whose solution is a  Brownian motion.

\bigskip

\noindent
{\bf III.  Heat equations, semigroups, and flows}

Let $\bar M$ be a compactification of $M$, i.e. a compact Hausdorff space.
 We assume $\bar M$ is first countable.  Let $h$ be a  continuous function 
 on $\bar M$. Consider the following  heat equation with initial boundary
 conditions:  
\begin{eqnarray}
	 \frac {\partial f} {\partial t}&=&{1 \over 2} \Delta f, 
\hskip 12pt  x\in M \\
	 f (x,0)&=&h(x), x\in\bar M\\
	f(x,t)&=&h(x), x\in \partial M
\end{eqnarray}

It is known that there is a unique minimal solution satisfying the first two
 equations on a stochastically complete manifold, the solution is in fact
 given by the semigroup associated with Brownian motion on the manifold
 applied to $h$. So  the above equation is not solvable in general. However
 with a condition imposed on the boundary of the  compactification,
 Donnelly-Li \cite{DO-LI84}  showed that the heat semigroup satisfies (4).
 Here is the condition and the theorem:

\noindent {\bf The ball convergence criterion:}
Let $\{x_n\}$ be a sequence in $M$ converging to a point $\bar x$ on the
 boundary,  then the geodesic balls $\B_r(x_n)$, centered at $x_n$ of
 radius $r$, converge to $\bar x$ as $n$ goes to infinity for each fixed $r$.  

An example of a manifold which satisfies the ball convergence criterion is
 $R^n$ with sphere at infinity, but not the  compactification of  a cylinder
 with a circle at infinity  added at each end. The one point compactification
 also satisfies the  ball convergence  criterion.

\begin{theorem}[H.Donnelly-P. Li] Let $M$ be a complete Riemannian manifold
 with Ricci curvature bounded from below. The over determined system (2)-(4)
 is solvable for any given continuous function  $h$ on $\bar M$,
if and only if the ball convergence criterion holds for $\bar M$.
\end{theorem}

Notice if the Brownian motion  starting  from $x$ converges to some
point on the boundary to which  $x$ converges, (2)-(4) is clearly solvable. 
See section 5 for details. We would also like to consider the opposite
 question: Do we get any information on the diffusion processes if we know
 the behaviour at infinity of the associated semigroups? This is true for
 many cases. In particular for the one point compactification, Schwartz has
 the following theorem\cite{SC89}, which provides a partial converse to
 \cite{EL82}:

\begin{theorem}[L. Schwartz]
 Let $\bar F_t$ be the standard extension of $F_t$ to
 $\bar M=M\cup \{\infty\}$, the one point compactification. Then the map
 $(t,x)\mapsto \bar F_t(x)$ is continuous from $R_+\times \bar M $ to
 $L^0(\Omega,P, \bar M)$, the space of measurable maps with topology induced
 from convergence in probability,   if and only if the semigroup $P_t$  has
 the $C_0$ property and the map $t\mapsto P_tf$ is continuous from $R_+$ to
 $C_0(M)$. 
\end{theorem}

\bigskip

\noindent{\bf{ Acknowledgement:}} 
This paper developed from my Warwick Msc. dissertation.
It is a pleasure to thank my supervisor professor D. Elworthy for pointing
out the  problem and for   discussions on the proofs and for
 critical reading of the paper. I am grateful to professor Zhankan Nie of Xi'an
 Jiaotong university and for the support of the Sino-British Friendship Scheme.
 I also benefited from the EC programme SCI-0062-C(EDB). I also like to thank
Prof. I. Chavel for helpful comments.

\section{Weak uniform cover and nonexplosion}

\begin{definition} \cite{ELbook} A stochastic dynamical system (1) is said to
 admit a uniform cover( radius $r>0$, bound $k$), if there are charts
 $\{\phi_n, U_n\}$ of diffeomorphisms from open sets $U_n$ of the manifold
 onto open sets $\phi_n(U_n)$ of $R^n$,  such that:
\begin{enumerate}
\item $B_{3r}\subset \phi_i(U_i)$, each $i$. ($B_\alpha$ denotes the open
 ball about $0$, radius $\alpha$).

\item The open sets $\{\phi_i^{-1}(B_r)\}$ cover the manifold.

\item  If $(\phi_i)_*(X)$ is given by:
$$(\phi_i)_*(X)(v)(e)=(D\phi_i)_{\phi_i^{-1}(v)}X(\phi_i^{-1}v)(e)$$
with  $(\phi_i)_*(A)$ similarly  defined, then both  $(\phi_i)_*(X)$ and 
 $\A(\phi_i)$   are bounded by $k$ on $B_{2r}$, here $\A$ is the generator
 for the dynamical system.
\end{enumerate}
\end{definition}

%\vskip 144pt

Let $\{e_j\}$ be an orthonormal basis for $R^m$, and define $X^j$ to be
$X(x)(e_j)$. Then  for Stratonovich equation (\ref{eq: basic}) on $M$, 
$$\A(\phi_i)=\sum_{ k,l=1}^m X^kX^l(\phi_i)+d\phi_i(A).$$
Let $M=R^n$. If we replace the Stratonovich differential in  (1) by
the  It\^ o differential, then
$$\A(\phi_i)(x)=\half\sum_{k,l=1}^m D^2\phi_i(x)(X^k(x),X^l(x)) +
 D\phi_i(x)(A(x)),$$
 which does not involve any of  the derivatives of $X^j$'s.

\begin{definition} 
   A diffusion process $\{F_t,\xi\}$ is said to have a weak uniform cover if 
there are pairs of connected open sets $\{U_n^0, U_n\}$,  and a non-increasing
 sequence   $\{\delta_n\}$ with $\delta_n>0$, such that:
\begin{enumerate}
\item 
$U_n^0\subset U_n$, and the open sets  $\{U_n^0\}$ cover the manifold.
 For $x\in U_n^0$ denote by  $\tau^n(x)$ the first exit time of $F_t(x)$ from 
	the open set $U_n$.  Assume $\tau^n<\xi$ unless $\tau^n=\infty$. 
\item  
  There exists $\{K_n\}_{n=1}^\infty$, a family of increasing open subsets of
   $M$ with  $\cup K_n=M$, such that each $U_n$ is contained in one of
 these sets and  intersects at most one boundary from
 $\{\partial K_m\}_{m=1}^\infty$.
\item 
 Let $x\in U_n^0$ and  $U_n\subset K_m$, then for $t<\delta_m$: 
\begin{equation}
 P\{\omega\colon \tau^{n}(x)< t\}\le Ct^2
\end{equation}
\item 
$\sum_{n=1}^\infty  \delta_n=\infty$. 
\end{enumerate}
\end{definition}

\vskip 144pt

Notice the introduction of $\{K_n\}$ is only for giving an order to the open
 sets $\{U_n\}$. This is  quite natural when looking at concrete manifolds.
 In a sense the condition  says  the geometry of the manifold under
 consideration changes slowly as far as  the diffusion process is concerned.
 In particular if  $\delta_n$ can be taken all  equal, we take $K_n=M$;
  e.g. when the number of open sets in the cover is  finite. On a Riemannian
 manifold the open sets are often taken as geodesic balls.

\noindent
\begin{lemma} Assume there is a uniform cover for the stochastic dynamic
 system(1), then the solution has a weak uniform cover with $\delta_n=1$,
 all $n$. 
\end{lemma}

\noindent
{\bf Proof:} This comes directly from lemma 5 , Page 127 in  \cite{ELbook}.

\bigskip
  
\noindent {\bf Remarks:}
\begin{enumerate}
\item
 Let $T$ be a stopping time, the inequality (5) gives the following  from
 the strong  Markov property of the process.: 

let $V\subset U_n^0$, and $V\subset K_m$, then when $t<\delta_m$:
\begin{equation}
P\{\tau^n(F_T(x))< t| F_T(x)\in V\}\leq Ct^2 ,
\end{equation}
since 
\begin{eqnarray*}
&P&\{\tau^n(F_T(x))<t, F_T(x)\in V\}\\
&=&\int_{V} P(\tau^n(y)<t)P_T(x,dy)
\leq Ct^2   P\{F_T(x)\in V\},
\end{eqnarray*}

here $P_T(x,dy)$ denotes the distribution of $F_T(x)$.
\item
Denote by $P_t^{U_n}$ the heat solution on $U_n$ with Dirichlet boundary condition, then (6) is equivalent to the following: when $x\in U_n^0$, 
\begin{equation}
1-P_t^{U_n}(1)(x)\leq Ct^2
\end{equation}
\item
The methods in this article work in infinite dimensions to give analogous
 results.
\end{enumerate}

\noindent
{\bf Exit times:\, }  Given such a cover, let $x\in U_n^0$. We define
 stopping times $\{T_k(x)\}$ as follows:  Let $T_0=0$. Let
 $T_1(x)=\inf\{t>0:  F_t(x,\omega)\not \in \overline {U_n}\}$
 be the first exit time of $F_t(x)$ from the set $U_n$. Then 
 $F_{T_1}(x,\omega)$ must be in one of the open sets $\{U_k^0\}$. Let
\begin{eqnarray*}
&\Omega_1^1&=\{\omega: F_{T_1}(x)(\omega)\in U_1^0, T_1(x,\omega)<\infty\}\\
&\Omega_k^1&=\{ \omega: F_{T_1}(x)\in U_k^0- \bigcup_{j=1}^{k-1} U_j^0, T_1<\infty\}
\end{eqnarray*}
Then $\{\Omega_k^1\}$ are disjoint sets such that
 $\cup \Omega_k^1=\{T_1<\infty\}$.
In general we only need to consider the nonempty sets of such. 
  Define further the following:
 Let $T_2=\infty$, if $T_1=\infty$. Otherwise if $\omega\in\Omega_k^1$, let: 
\begin{equation}
T_2(x,\omega)=T_1(x,\omega)+\tau^k(F_{T_1}(x,\omega))
\end{equation}

In a similar way, the whole sequence of stopping times $\{T_j(x)\}$ and sets 
$\{\Omega_k^j\}_{k=1}^\infty$ are defined for $j=3, 4, \dots$. Clearly
 $\Omega_k^j$ so defined is measurable with respect to the sub-algebra
 $\F_{T_j}$.

\begin{lemma}
Given a weak uniform cover as above. Let $x\in U_n^0$ and $U_n\subset K_m$. Let 
 $t< \delta_{m+k}$. Then
\begin{equation}
P\{\omega: T_k(x,\omega)-T_{k-1}(x,\omega)< t, T_{k-1}<\infty\}\leq Ct^2
 \end{equation}  
\end{lemma}

\noindent {\bf Proof:} Notice for such  $x$, $F_{T_{k}}(x)\in K_{m+k-1}$. Therefore for $t< \delta_{m+k}$ we have: 
\begin{eqnarray*}
&P&\{\omega: T_k(x,\omega)-T_{k-1}(x,\omega)< t, T_{k-1}<\infty\}\\
&=& \sum_{j=1}^\infty P(\{\omega: T_k(x,\omega)-T_{k-1}(x,\omega)< t\}
\cap \Omega^{k-1}_j)\\
&=& \sum_{j=1}^\infty P\{\tau^j(F_{T_{k-1}}(x))<t, \Omega^{k-1}_j\}\\
&\leq& Ct^2\sum_{j=1}^\infty P(\Omega^{k-1}_j)\leq Ct^2, 
\end{eqnarray*}           
as in remark 1. Here $\chi_A$ is the characteristic function for a measurable
 set $A$, and  $E$ denotes taking expectation. 
\hfill\rule{3mm}{3mm}

\begin{lemma} If $\sum_n t_n=\infty$, $t_n>0$  non-increasing. Then there is a
non-increasing sequence $\{s_n\}$, such that $0<s_n\le t_n$:

(i) $\hskip 12pt \sum s_n=\infty$

(ii) $\hskip 12pt \sum s_n^2 <\infty$
\end{lemma}
{\bf Proof:}
 Assume $t_n\le 1$, all $n$. Group the sequence $\{t_n\}$ in the following
 way: 
$$t_1;\ t_2,\ \dots ,\ t_{k_2};\ t_{k_2+1},\ \dots , \ t_{k_3};
 \ t_{k_3+1}\ \dots$$
Such that $1\le t_2+\dots+t_{k_2}\le 2$,\ $1\le t_{k_i+1} + 
\dots +t_{k_{i+1}}\le 2$, \ $i\ge 2$.
\ Let $s_1=t_1$, $s_2=\frac {t_2} 2, \dots $, $s_{k_2}=\frac {t_{k_2}} 2$,
 $s_{k_2+1}=\frac {t_{k_2+1}} 3,\dots$, $s_{k_3}=\frac {t_{k_3}} 3$,
 $s_{k_3+1}=\frac {t_{k_3+1}} 4, \dots$.  \ Clearly the $s_n$'s so defined
 satisfy the requirements.
\hfill\rule{3mm}{3mm}

So without losing generality,  we may assume from now on that the constants 
$\{\delta_n\}$ for a weak uniform cover fulfill the two conditions in the 
above lemma. With these established, we can now state the nonexplosion result.
Th proof is  analogous to the  argument of theorem  6 on Page 129 
 in \cite{ELbook}.

\begin{theorem} If the solution $F_t(x)$ of the equation $(1)$ has a weak
 uniform cover,  then it is  complete(nonexplosion).
\end{theorem}
{\bf Proof: \ }
Let $x\in K_n$, $t>0$, $0<\epsilon< 1$. Pick  a number $p$ (may be depending
 on $\epsilon$ and $n$), such that $\sum_{i=n+1}^{n+p}\epsilon\delta_i>t$,
 which is possible from:  $\sum_{i=1}^\infty \delta_i=\infty$.
So
\begin{eqnarray*}
 P\{\xi(x)< t\}&\leq&     P\{ T_p(x)< t, T_{p-1}<\infty\}\\
&=&    P\{\sum_{k=1}^p (T_k(x)-T_{k-1}(x))<t, T_{p-1}<\infty\}\\
&\leq& \sum_{k=1}^p P\{ T_k(x)-T_{k-1}(x)< \epsilon \delta_{n+k}, T_{k-1}<\infty\} \\
&\leq& C\epsilon^2\sum_{k=n}^{n+p} \delta_k^2   
\leq C\epsilon^2\sum_{k=1}^{\infty} \delta_k^2.
\end{eqnarray*}
Let $\epsilon\to 0$, we get:  $P\{ \xi(x)<t\}=0$. 
\hfill\rule{3mm}{3mm}

\bigskip

\noindent{\bf Remark:  }   
The argument in the above proof is valid if the definition 
of a weak uniform cover is changed slightly, i.e.  replacing the constant $C$ 
by $C_n$ (with some slow growth condition) but keep all $\delta_n$ equal.

\bigskip

To understand more of the weak uniform cover, we look into an example:

\noindent
{\bf Example 1 :}  Let $\{U_n\}$ be a family of relatively compact open
 (proper) subsets of $M$ such that $U_n\subset U_{n+1}$ and
 $\cup_{n=1}^\infty U_n=M$. Assume there is a 
sequence of numbers $\{\delta_n\} $ with $\sum_n \delta_n=\infty$, such that the 
following inequality holds  when $t<\delta_{n-1}$ and $x\in U_{n-1}$:
$P\{\tau^{U_n}(x)<t\}\leq ct^2$. 
Then the diffusion concerned does not explode by taking
 $\{U_{n+1}-\overline{U_{n-1}},U_n-\overline{U_{n-1}}\}$ to be a weak uniform cover and $K_n=\bar U_n$.

\section{Boundary behaviour of diffusion processes}

To consider the boundary behaviour of  diffusion processes, we introduce the 
following concept:

\begin{definition} A weak uniform  cover $\{U_n^0, U_n\}$ is said to be
 regular (at infinity for $\bar M$), if the following holds: let $\{x_j\}$ be
 a sequence in $M$ converging to $\bar x\in \partial M$, and
$x_j\in U_{n_j}^0\in\{U_n^0\}_{n=1}^\infty$, then the corresponding open sets 
$\{U_{n_j}\}_{j=1}^\infty  \subset \{U_n\}_1^\infty$ converges to $\bar x$ as
 well.   A regular uniform cover can be  defined in a similar way. 
\end{definition}
 
For a point $x$ in $M$, there are a succession of related open sets 
$\{W_x^p\}_{p=1}^\infty$,  which are defined as follows:
Let $W_x^1$ be the union of all open sets from $\{U_n\}$ such that  $U_n^0$
contains $x$, and $W_x^2$ be the union of all open sets from 
$\{U_n\}$ such that $U_n^0$ intersects one of the open sets $U_{n_j}^0$
 defining $W_x^1$. The $\{W_{x}^{p}\}$ are defined similarly. These sets are
 well defined and in fact form an increasing sequence.

\begin{lemma} Assume $F_t$ has a regular weak uniform cover. Let $\{x_n\}$ be
 a  sequence in $M$  which converges to a point $\bar x \in \partial M$.
 Then $W_{x_n}^p$ converges to  $\bar x$ as well for each fixed $p$.
\end{lemma}

\noindent
{\bf Proof:} We only need to prove the following:
let  $\{z_k\}$ be a sequence  $z_k\in W_{x_k}^p$, then $z_k \to \bar x$, as
 $n\to \infty$.  First let $p=2$. 

By definition, for each $x_k$, $z_k$, there are open
 sets $U_{n_k}^0$ and $U_{m_k}^0$ such that $x_k\in U_{n_k}^0$, 
  $z_k\in U_{m_k}^0$ and $U^0_{n_k}\cap U^0_{m_k}\not = \emptyset$. Furthermore $U_{n_k}\to \bar x$ as $k\to \infty$. 
Let $\{y_k\}$ be a sequence of points with $y_k\in U_{n_k}^0\cap U_{m_k}^0$.  But $y_k\to \bar x$ since $U_{n_k}$ does.  So $U_{m_k}\to \bar x$ again 
from the definition of a regular weak uniform cover. Therefore $z_k$ converges to $\bar x$ as $k\to \infty$, which is what we want. 
The rest can be proved by induction.
\hfill\rule{3mm}{3mm}

\begin{theorem}
 If the diffusion $F_t$ admits a regular weak uniform cover for $\bar M$, with
 $\delta_n=\delta$, all $n$, then the map $F_t(-):M\to M$ can be extended to
 the compactification $\bar M$ continuously in probability with the
 restriction to the boundary to be the identity map,  uniformly in $t$ in
 finite intervals. (We will say $F_t$ extends.)
\end{theorem}

\noindent{\bf Proof:} 
  Take $\bar x\in \partial M$ and a sequence $\{x_n\}$ in $M$ converging to
 $\bar x$.  Let $U$ be a neighbourhood of $\bar x$ in $\bar M$. We want to
 prove for each $t$: 

$$\lim_{n\to\infty}P\{\omega:F_s(x_n,\omega)\not\in U,
 \hskip 5pt \hbox{for some}\hskip 3pt s<t \}=  0.$$
 Since $x_n$ converges to $\bar x$, 
	there is a number $N(p)$ for each $p$, such that if $n>N=N(p)$,
 $W_{x_n}^p\subset U$. Let  $t>0$, choose $p$ such that ${2t \over p}<\delta$.
	For a number $n>N(p)$ fixed,   we have:
\begin{eqnarray*}
 &P&\{\omega: F_s(x_n,\omega)\not \in U,
\hskip 5pt \hbox{for some}\hskip 3pt s< t\}\\
 &\subset& P\{\omega: F_s(x_n,\omega)\not\in W_{x_n}^p,
 \hskip 5pt \hbox{for some}\hskip 3pt s< t\}\\
	&\leq& P\{\omega:T_p(x_n)(\omega)< t,T_{p-1}(x_n)<\infty\}\\
	&\leq& \sum_{k=1}^{p} P\{T_k(x_n)-T_{k-1}(x_n)<  {t \over p} , 
		T_{k-1}(x_n)<\infty\}\\
	&\leq& {Ct^2 \over p}	\end{eqnarray*}
Here $C$ is the constant in the definition of the weak uniform cover.
 Let $p$ go  to infinity to complete the proof.
\hfill\rule{3mm}{3mm}

\bigskip

\noindent {\bf Remark:}
If $\delta_n$ can be taken all equal, theorem 2.4, theorem 3.2  hold if (5)
 is relaxed to:         $$P\{\tau^n(x)<t\}\le f(t),$$
for some nonnegative function $f$ satisfying $\lim_{t\to 0} {f(t) \over t}=0$. 

\bigskip

 The weak uniform cover condition is not very easy to apply for a general
 diffusion on a manifold. For a  Brownian motion in a Riemannian manifold it is
 relatively simple. We often start with a uniform cover, ref. lemma 2.1.

\bigskip

\noindent
{\bf Example 2:\ }
Let $M=R^n$ with the one point compactification. Consider the following s.d.e.:
$$ (It\widehat o) \hskip 12pt  dx_t=X(x_t)dB_t+A(x_t)dt$$ 
Then if both $X$ and $A$ have linear growth, the solution has the $C_0$
 property.

{\bf Proof:}
 There is a well known uniform cover for this system. See \cite{Clark73},
 or \cite{ELbook}. 
 A slight change gives us the following regular uniform cover:

Take a countable set of points $\{p_n\}_{n\ge 0}\subset M$ such that the open sets 
 $U_n^0=\{z:|z-p_n|<{ |p_n| \over 3}\}, n=1,2, \dots$  and  $U_0^0=\{ z: |z-p_0|<2\}$
cover $R^n$. Let $U_0=\{ z: |z-p_0|<6\}$; and $U_n=\{z:|z-p_n|<{|p_n|\over 2}\}$, for $n\not = 0$.  Let  $\phi_n$ be the chart map on $U_n$:
$$\phi _n(z)={z-p_n\over |p_n|}.$$
This certainly defines a uniform cover(for details see Example 3 below). Furthermore if $z_n\to \infty$ and
$z_n\in U_n^0$, then any $y\in U_n$ satisfies the following:
 $$|y|>{|p_n| \over 2} >{1\over 3}|z_n| \to \infty,$$
since $|p_n|\ge {{3 |z_n|} \over 4}$.
Thus we have a regular uniform cover which gives the required $C_0$ property.

\bigskip

\noindent {\bf Example 3:}  
Let $M=R^n$, compactified with the sphere at infinity: 
 $\bar M=R^n\cup S^{n-1}$. Consider the same s.d.e as in the example above. 
Suppose both $X$ and $A$ have sublinear growth of power $\alpha<1$: 
\begin{eqnarray*}
|X(x)|&\le&   K(|x|^\alpha+1)\\
 |A(x)|& \le& K(|x|^\alpha+1)
\end{eqnarray*}
for a constant $K$. Then there is no explosion. Moreover the solution 
 $F_t$ extends.

\noindent {\bf Proof:} 
 The proof is as in example 2, we only need to construct a regular
 uniform cover for the s.d.e.:

Take points  $p_0, p_1, p_2, \dots$ in $R^n$ (with $|p_0|=1$, $|p_n|>1$) ,  such that  the open sets $\{U_n^0\}$ defined  by : 
 $ U_0=\{ z: |z-p_0|<2\}$, $U_i^0=\{z:|z-p_i|<{|p_i|^\alpha \over 6}\}$ cover $R^n$. 

 Let $U_0=\{ z: |z-p_0|<6\}$, $U_i=\{z:|z-p_i|<{|p_i|^\alpha\over 2}\}$, and
let   $\phi_i$ be the chart map from $U_i$ to $R^n$:
$$\phi _i(z)={6(z-p_i)\over |p_i|^\alpha}.$$

Then $\{\phi_i, U_i\}$ is a uniform cover for the stochastic dynamical system. 
In fact, for $i\not = 0$, and $y\in B_3\subset R^n$: 
$$ |(\phi_i)_\ast (X)(y)|\le {K(1+|\phi_i^{-1}(y)|^\alpha) \over |p_i|^\alpha}
\leq {K \over |p_i|^\alpha} (1+2|p_i|^\alpha)<18K$$

Similarly $|(\phi_i)_\ast (A)(y)|\le 18K$, and $D^2\phi_i=0$.

Next we show this cover is regular. Take a sequence  $x_k$ converging to  $\bar x$
in $\partial R^n$. Assume $x_k\in U_k^0$. Let $z_k\in U_k$. We want to prove $\{z_k\}$ converge to $\bar x$.  First the norm of $z_k$ converges to infinity as $k \to \infty$, since $|p_k|>{ {2|x_k|} \over 3}$ and $|z_k|>|p_k|-{1\over 2} |p_k|^\alpha$. 

\vskip 144pt

Let $\theta$ be the biggest angle between points in $U_k$, then 

$$  \hbox{tan}{\theta\over 2} \leq \sup_{z\in U_n} { |z-p_n| \over |p_n|}
\le {|p_n|^\alpha \over 2|p_n|}
\le {|p_n|^{\alpha -1}\over 2} \to 0.         $$
Thus $\{U_n, \phi_n\}$ is a uniform cover satisfying the convergence criterion
 for the sphere compactification. The required result holds from the theorem. 

\bigskip

This result is sharp in the sense there is a s.d.e. with coefficients having
linear growth but the solution to it does not extend to the sphere at
 infinity to be identity:

\noindent {\bf Example 4:} 
 Let $B$ be a one dimensional Brownian motion. Consider the 
following s.d.e on the complex plane $\cal C$:
		$$dx_t=ix_tdB_t$$
The solution starting from $x$ is in fact $xe^{iB_t+{t\over 2}}$, which does
 not continuously extend to be the identity on the sphere at infinity.

\section{Boundary behaviour continued}

A diffusion process is a $C_0$ diffusion if its semigroup has the $C_0$
 property.  This is equivalent to the following\cite{AZ74}: let $K$ be a
 compact set, and $T_K(x)$ the first entrance time to $K$ of the diffusion
 starting from $x$, then  $\lim_{x\to\infty} P\{T_K(x)<t\}=0$ for each $t>0$,
 and each compact set $K$. 

The following theorem follows from theorem 3.2 when $\delta_n$ in the
 definition of  weak uniform cover can be taken all equal:

\begin{theorem} Let $\bar M$ be the one point compactification. Then if the 
diffusion process $F_t(x)$ admits a regular weak uniform cover, it is a 
$C_0$ diffusion.
\end{theorem}

{\bf Proof:}
Let $K$ be a compact set with $K\subset K_j$; here $\{K_j\}$ is as in
 definition 2.2.
Let $\epsilon>0$, $t>0$, then there is a number $N=N(\epsilon, t)$ such that:

$$\delta_{j+2} +\delta_{j+4}+\dots+ \delta_{j+2N-2}>{t\over \epsilon}$$

Take $x\not \in K_{j+2N}$.
Assume $x\in K_m$, some $m>j+2N$. Let $T_0$ be the first entrance time of 
$F_{t}(x)$ to $K_{j+2N-1}$,   $T_1$ be the first entrance time  of 
$F_{t}(x)$ to $K_{j+2N-3}$ after $T_0$, (if $T_0<\infty$), and so on. 
 But $P\{T_i<t, T_{i-1}<\infty\}\le Ct^2$ 
for $t<\delta_{j+2N-2i}, i>0$, since any open sets from the cover  intersects at most one boundary of sets from $\{K_n\}$. Thus
\begin{eqnarray*}
P\{T_K(x)<t\}&\leq&   P\{\sum_1^{N-1} T_{i}(x)<t, T_{N-2}<\infty\}   \\
&\leq& \sum_{i=1}^{N-1} P\{T_i(x)<\epsilon \delta_{j+2N-2i}, T_{i-1}(x)<\infty\}\\
&\le& C\epsilon^2\sum_1^{N-1}\delta_{j+2N-2i}^2 \, 
\le C\epsilon^2\sum_1^\infty\delta_{j}^2
\end{eqnarray*}
The proof is complete by letting $\epsilon \to 0$.
\hfill\rule{3mm}{3mm} 

\noindent
{\bf Example 5:\,} \label{ex: nonexplosion fo Ricci}
 Let $M$ be a complete Riemannian manifold, $p$ a fixed
 point in $M$. Denote by $\rho(x)$ the distance between $x$ and $p$, 
 ${\cal B}_r(x)$ the geodesic ball centered at $x$ of  radius  $r$, and
 Ricci$(x)$ the Ricci curvature at $x$.

\noindent {\bf Assumption A:}  
\begin{equation}
\int_1^\infty {1\over \sqrt{K(r)}} dr=\infty   
 \end{equation}

\noindent
Here $K$ is defined as follows:

$$K(r)=-\{\inf_{B_r(p)} \hbox{Ricci}(x)\wedge 0\}$$

Let $X_t(x)$ be a Brownian motion on $M$ with $X_0(x)=x$. Consider the first
 exit  time of $X_t(x)$ from $B_1(x)$:

 $$T=\inf_t\{t\geq 0: \rho(x ,X_t(x))=1\}.$$ 

Then we have the following  estimate on $T$ from \cite{Hsu89}:

If $L(x)>\sqrt{K(\rho(x)+1)}$, then
$$P\{T(x)\leq {c_1 \over L(x)} \}\leq \exp^{-c_2L(x)}$$
for all $x\in M$. Here $c_1, c_2$ are positive constants independent of $L$. 

\noindent
This can be rephrased into the form we are familiar with:
 when $0\le t<{c_1 \over \sqrt{K(\rho(x)+1)}}$,
$$P\{T(x)\leq t\}\leq \exp^{-{c_1c_2   \over t}}$$
 But  $\lim_{t\to 0} {e^{-{c_1c_2\over t}} \over t^2}=0$. So there is a 
 $\delta_0>0$, such that:
$e^{-{c_1 c_2 \over t}} \leq t^2$, when $t<\delta_0$.
Thus:
 
\noindent
{\bf Estimation on exit times:  \,} 
when $t<{c_1 \over \sqrt{K(\rho(x)+1)}}\wedge \delta_0 $,
\begin{equation}
P\{T(x)<t\}\leq t^2.
\end{equation}

\noindent
Let $\delta_n={1\over \sqrt{K(3n+1)}}\wedge \delta_0$, then we also have the
 following :
\begin{equation}
\sum_1^\infty \delta_n      
\ge   \sum_1^\infty{1\over \sqrt{K(3n+1)}}
\ge \int_1^\infty {1\over \sqrt{K(3r)}}\, dr=\infty
\end{equation}

with this we may proceed to prove the following from Hsu:

\noindent {\bf Corollary:} [Hsu] \,\,
A complete Riemannian manifold $M$ with Ricci curvature satisfying assumption A
 is  stochastically complete and has the $C_0$ property.

\noindent {\bf Proof:} 
There is a regular weak uniform cover as follows:

 First take any $p\in M$, and let $K_n=\overline{B_{3n}(p)}$. 
Take points $p_i$ such that $U_i^0=B_1(p_i)$ covers the manifold. Let $U_i=B_2(p_i)$. Then  $\{U_i^0,U_i\}$ is a regular weak uniform cover for $M\cup \Delta$.  

\bigskip

\noindent{\bf Remark:} Grigory\'an has the following volume growth test on nonexplosion.  The Brownian motion does not explode on a manifold if there is a function $f$
on $M$ satisfying: 
$$\int^\infty {r\over \hbox{Ln}(\hbox{Vol}(B_R))}\, dr=\infty.$$
Here $\hbox{Vol}(B_R)$ denotes the volume of a geodesic ball centered at a point $p$ in  $M$. 
This result is stronger than the corollary obtained above by the following 
comparison theorem on a $n$ dimensional manifold: let $\omega_{n-1}$ denote the 
volume of  the $n-1$ sphere of radius $1$,
 
$$\hbox{Vol}(B_R) \le \omega_{n-1} \int_0^R\{ \sqrt{{(n-1) \over K(R)}} 
\hbox{Sinh}(\sqrt{{K(R)\over (n-1)}} r)\}^{(n-1)}\, dr.$$
Notice $K(R)$ is  positive when $R$ is sufficiently big provided the  Ricci
 curvature is not nonnegative everywhere. So Grigory\'an's result is stronger
 than the one obtained above.

The definition of weak uniform cover is especially suitable for the one point 
compactification. For general compactifications the following definition 
explores more of the geometry of the manifolds and gives a better result:

\begin{definition}
Let $\bar M$ be a compactification of $M$, $\bar x \in \partial M$. A
 diffusion process $F_t$ is said to have a uniform cover at the  point
 $\bar x$, if there is a sequence  $A_n$ of open neighbourhoods of $\bar x$
 in $\bar M$ and positive numbers $\delta_n$  and a constant $c>0$, such that:
\begin{enumerate}
\item 
The sequence of $A_n$ is strictly decreasing, with $\cap A_n=\bar x$, and 
$A_n\supset \partial A_{n+1}$. 
\item  
The sequence of numbers $\delta_n$ is non-increasing  with
 $\sum \delta_n=\infty$ and $\sum\delta^2_n<\infty$.
\item When $t<\delta_{n}$, and $x\in A_{n}-A_{n+1}$,
$$P\{\tau^{A_{n-1}}(x)<t\}\leq ct^2$$
Here $\tau^{A_n}(x)$ denotes the first exit time of $F_t(x)$ from the set $A_n$.
\end{enumerate}
\end{definition}
 
%\vskip 144 pt

\begin{proposition}
If there is a uniform cover for $\bar x \in \partial M$, then $F_t(x)$
 converges to $\bar x$ continuously in probability,
 uniformly in $t$ in finite interval,  as $x \to\bar x$.
\end{proposition}

Proof: 
The existence of $\{A_n\}$ will ensure $F_{\tau^{A_{n}}}(x)\subset A_{n-1}$, 
which allows us to apply a similar argument as in the case of the one point 
compactification. Here we denote by $\tau^A$ the first exit time of the process 
$F_t(x)$ from a set $A$.

Let $U$ be a neighbourhood of $\bar x$. For this  $U$, by compactness of $\bar M$, there is a number $m$ such that $A_m\subset U$, since $\cap_{k=1}^\infty A_k=\bar x$.
Let $0<\epsilon<1$, $\bar \epsilon=({\epsilon \over c\sum_k \delta_k^2})^{1\over 2}$.
we may assume $\bar \epsilon<1$.  Choose $p=p(\epsilon)>0$ such that: 

$$\delta_{m}+\delta_{m+1}+\dots +\delta_{m+p-1}>{t \over\bar \epsilon}$$

Let $x\in A_{m+p+2}$. Denote by $T_0(x)$ the first exit time of $F_t(x)$ from
 $A_{m+p+1}$, $T_1(x)$ the first exit time of $F_{T_0}(x)$ from $A_{m+p}$ where 
defined.  Similarly $T_i, i>1$ are defined.

Notice  if $T_i(x)<\infty$, then $F_{T_i(x)}\in A_{m+p-i}-A_{m+p+1-i}$, for $i=0, 1, 2,\dots $. Thus for $i>0$ there is the following inequality from the definition and the Markov property:
$$P\{T_i(x)<\bar\epsilon \delta_{m+p-i}\}\leq c\bar\epsilon^2 \delta^2_{m+p-i}$$

Therefore we have:
\begin{eqnarray*}
P\{\tau^U(x)<t\}    &\leq&  P\{\tau^{A_m}(x)<t\}     \\
&\leq&P\{T_{p}+\dots +T_1<t, T_{p-1}<\infty\}\\
&\leq& \sum_{i=1}^{p} P\{T_i<\bar \epsilon\delta_{m+p-i}, T_{i-1}<\infty\}\\
&\leq& c\bar\epsilon^2\sum_{i=1}^{p}\delta_{m+p-i}^2<\epsilon     
\end{eqnarray*}
This finishes the proof.   \hfill\rule{3mm}{3mm}

\section{Property at infinity of semigroups}
Recall a semigroup is said to have  the $C_0$ property, if it sends 
$C_0(M)$, the space of continuous functions on $M$ vanishing at infinity, to itself. 
Let $\bar M$ be a compactification of $M$. Denote by $\Delta$ the  point at infinity
for the one point compactification.  Corresponding to the $C_0$ property of
 semigroups we consider the following $C_*$ property for $\bar M$: 

\begin{definition}
A semigroup $P_t$ is said to have the $C_*\hskip 4pt$ property for $\bar M$, 
if for each continuous function $f$ on $\bar M$, the following holds: let
 $\{x_n\}$ be a sequence converging to $\bar x$ in $\partial M$, then
\begin{equation}
	\lim_{n\to \infty} P_tf(x_n)=f(\bar x), 
\end{equation}
\end{definition}

To justify the definition, we notice if $\bar M$ is the one point
 compactification, condition $C_*\hskip 4pt$ will imply the $C_0\hskip 4pt$
 property of the semigroup.  On the other hand if $P_t$ has the $C_0$ property,
 it has the $C_*$ property for  $M\cup \Delta$ assuming nonexplosion. This is
 observed by subtracting a constant function from a continuous function $f$ on
 $M\cup\Delta$: Let $g(x)=f(x)-f(\Delta)$, then $g\in C_0(M)$.
 So $P_tg(x)=P_tf(x)-f(\Delta)$. Thus 
$$\lim_{n \to \infty} P_tf(x_n)
=\lim_{n \to \infty} P_tg(x_n) +f(\Delta)=f(\Delta),$$ 
if $\lim_{n \to \infty} x_n=\Delta$. 

In fact the $C_*$ property holds for the one point compactification if and
 only if there is no explosion and the $C_0$ property holds. These properties
 are often   possessed by processes, e.g. a Brownian motion on a Riemannian
  manifold with  Ricci curvature which satisfies (10) has this property.

Before proving this claim, we observe first that:

\begin{lemma}
 If $P_t$ has the $C_*$ property for any compactification $\bar M$, 
it must have the $C_*$ property for the one point compactification. 
\end{lemma}

\noindent {\bf Proof:\ }
Let $f\in C(M\cup \Delta)$. Define a map $\beta$ from $\bar M$ to
 $M\cup \Delta$:  $\beta(x)=x $ on  the interior of $M$, and $\beta(x)=\Delta$,
 if $x$ belongs   to the boundary. Then $\beta$ is a continuous map from
 $\bar M$ to  $M\cup \Delta$, since for any compact set $K$, the inverse set
 $\bar M-K=\beta^{-1}(M\cup\Delta-K)$  is open in $\bar M$.

Let $g$ be the composition map of $f$ with $\beta$:
$g=f\circ \beta: \bar M\to R$.  Thus $g(x)|M=f(x)|M$,  and
 $g(x)|_{\partial M}=f(\Delta)$. So for a sequence $\{x_n\}$ converging to
 $\bar x\in \partial M$,
$\lim_n P_tf(x_n)=\lim_n P_tg(x_n)=g(\bar x)=f(\Delta)$.               
\hfill\rule{3mm}{3mm}

We are ready to prove the following theorem:
\begin{theorem}
If a semigroup $P_t$ has the $C_* $ property, the associated diffusion process $F_t$ is complete.
\end{theorem}
Proof:  We may assume the compactification under consideration is the one point
compactification from the the lemma above. Take $f\equiv 1$, $P_tf(x)=P\{t<\xi(x)\}$.
But $P\{t<\xi(x)\}\to 1$ as $x\to \Delta$ from the assumption. 
More precisely for any $\epsilon>0$, there is a compact set $K_\epsilon$ such that
if $x\not\in K_\epsilon$, $P\{t<\xi(x)\}>1-\epsilon$.

Let $K$ be a compact set containing $K_\epsilon$. Denote by $\tau $ the first
exit time of $F_t(x)$ from $K$. So $F_\tau(x)\not \in K_\epsilon$ on the set 
$\{\tau<\infty\}$. Thus:
\begin{eqnarray*}
P\{t<\xi(x)\}
 &\geq& P\{\tau<\infty, t<\xi(F_\tau(x))\}+P\{\tau=\infty\}\\
 &=& E\{\chi_{\tau<\infty} E\{\chi_{t<\xi(F_\tau(x))}|\F_\tau\}\}+  P\{\tau=\infty\}
\end{eqnarray*}
Here $\chi_A$ denotes the characteristic function of set $A$. Applying the strong Markov property of the diffusion we have:
$$E\{\chi_{\tau<\infty} E\{\chi_{t<\xi(F_\tau(x))}|\F_\tau\}\}
=E\{\chi_{\tau<\infty} E\{t<\xi(y)|F_\tau=y\}\}.$$
However $$E\{\chi_{t<\xi(y)}|F_\tau=y\}> 1-\epsilon,$$
So
\begin{eqnarray*}
P\{t<\xi(x)\}&\ge& P\{\tau=\infty\}+E\{\chi_{\tau<\infty}(1-\epsilon)\}  \\
&=& 1-\epsilon P\{\tau<\infty\}
\end{eqnarray*}
Therefore  $P\{t<\xi(x)\}=1$, since $\epsilon$ is arbitrary.
\hfill\rule{3mm}{3mm}

\bigskip

In the following we examine the  the relation between the behaviour at $\infty$
of  diffusion processes  and the diffusion semigroups.

\begin{theorem}
The semigroup $P_t$ has the $C_*$ property if and only if the diffusion process
$F_t$ is complete and can be extended to $\bar M$ continuously in probability
with $F_t(x)|_{\partial M}=x$.
\end{theorem}
Proof:
    Assume $F_t$ is complete and extends. 
  Take a point $\bar x \in \partial M$, and sequence $\{x_n\}$ converging to
 $\bar x$. Thus 
$$\lim_{n\to \infty} P_tf(x_n)=\lim_{n\to \infty} Ef(F_t(x_n))=Ef(x)=f(x)$$
 for any continuous function on $\bar M$, by the dominated convergence theorem.
	
 On the other hand $P_t$ does not have the $C_*$ property if the  assumption
 above is not true. In fact let $x_n$ be a sequence converging to $\bar x$,
 such that for some neighbourhood $U$ of $\bar x$, and a number $\delta>0$:
	$$\liminf_{n\to \infty} P\{F_t(x_n)\not \in U\}=\delta$$
 There is therefore a subsequence  $\{x_{n_i}\}$ such that:
	$$\lim_{i\to \infty} P\{F_t(x_{n_i})\not \in U\}=\delta.$$
 Thus there exists $N>0$, such that if $i>N$:
	$$P\{F_t(x_{n_i}) \in \bar M-U\}>{\delta \over 2}.$$
But since $\bar M$ is a compact Hausdorff space, there is a continuous 
function $f$ from $M$ to $[0,1]$ such that $f|_{\bar M-U}=1$, and $f(x)|_G=0$,
 for any  open set $G$ in $U$. Therefore
	
\begin{eqnarray*}
P_tf(x_n)&=&Ef(F_t(x_n))  \\
&\geq& \int_{\{\omega: F_t(x_n)\in \bar M-U\}}
	 f(F_t(x_n))\, P(d\omega)\\
&=&P\{F_t(x_n)\in \bar M-U\} >  {\delta \over 2}  
\end{eqnarray*}
So $\lim P_tf(x_n)\not =f(x)=0.$ 
\hfill\rule{3mm}{3mm}
	
\begin{corollary}
Assume the diffusion process $F_t$ admits a weak uniform cover regular for
 $M\cup \Delta$, then its diffusion semigroup $P_t$ has the $C_*$ property 
for $M\cup\Delta$. The same is true 
 for a general compactification if all $\delta_n$ in the weak uniform cover 
can be taken equal.
\end{corollary}

\noindent {\bf Example 7:}\cite{DO-LI84}\,
Let $M$ be a complete connected Riemannian manifold with Ricci curvature
 bounded from below. Let  $\bar M$ be  a compactification such that the ball
 convergence criterion holds (ref. section 1).  In particular the over
 determined equation (2)-(4) is  solvable  for any continuous function $f$
 on $\bar M$ if  the ball convergence  criterion holds.

\noindent {\bf Proof:} 
We keep the notation of example 5 here.  Let $K=-\{\inf_x Ricci(x)\wedge 0\}$,
 $\delta={c_1\over k}$, where $c_1$ is the constant in example 5. Let
 $p\in M$ be a fixed point, and $K_n=\overline{B_{3n}(p)}$ be compact sets in
 $M$.  Take points $\{p_i\}$ in $M$ such that $\{B_1(p_i)\}$ cover the
 manifold.  Then $\{B_1(p_i), B_2(p_i)\}$ is a weak uniform cover from (11)
 and (12). Moreover  this is a regular cover if the ball convergence
 criterion holds for the compactification.

%\bibliographystyle{plain}
%\bibliography{paper,paper2}

\begin{thebibliography}{10}

\bibitem{AZ74}
R.~Azencott.
\newblock Behaviour of diffusion semigroups at infinity.
\newblock {\em Bull. Soc. Math, France}, 102:193--240, 1974.

\bibitem{BA86}
D.~Bakry.
\newblock Un crit\'ere de non-explosion pour certaines diffusions sur une
  var\'iet\'e {R}iemannienne compl\'ete.
\newblock {\em C. R. Acad. Sc. Paris}, t. 303, S\'erie I(1):23--26, (1986).

\bibitem{CA-EL83}
A.~Carverhill and D.~Elworthy.
\newblock Flows of stochastic dynamical systems: the functional analytic
  approach.
\newblock {\em Z. Wahrscheinlichkeitstheorie Verw. Gebiete}, 65:245--267,
  (1983).

\bibitem{Clark73}
J.M.C. Clark.
\newblock An introduction to stochastic differential equations on manifolds.
\newblock In D.Q. Mayne;~R.W. Brockett, editor, {\em Geometric methods in
  systems theory}, pages 131--149, (1973).

\bibitem{DA89}
E.B. Davies.
\newblock Heat kernel bounds, conservation of probability and the {\sc f}eller
  property.
\newblock {\em preprint}, (1989).

\bibitem{Dodz}
J.~Dodziuk.
\newblock Maximal principle for parabolic inequalities and the heat flow on
  open manifolds.
\newblock preprint.

\bibitem{DO-LI84}
H.~Donnelly and P.~Li.
\newblock Heat equation and compactifications of complete {\sc r}iemannian
  manifold.
\newblock {\em Duke mathematical journal}, 51, No. 3:667--673, 1984.

\bibitem{ELflour}
K.~D. Elworthy.
\newblock Geometric aspects of diffusions on manifolds.
\newblock In P.~L. Hennequin, editor, {\em Ecole d'Et\'e de Probabilit\'es de
  Saint-Flour XV-XVII, 1985, 1987. Lecture Notes in Mathematics}, volume 1362,
  pages 276--425. Springer-Verlag, 1988.

\bibitem{ELbook}
K.D. Elworthy.
\newblock {\em Stochastic {D}ifferential {E}quations on {M}anifolds}.
\newblock Lecture Notes Series 70, Cambridge University Press, 1982.

\bibitem{EL82}
K.D. Elworthy.
\newblock Stochastic flows and the $c_0$ property.
\newblock {\em Stochastics}, 0:1--6, (1982).

\bibitem{GA59}
M.P. Gaffney.
\newblock The conservation property of the heat equation on {R}iemannian
  manifolds.
\newblock {\em Comm. Pure. appl. math.}, 12, 1959.

\bibitem{GR87}
A.~Grigoryan.
\newblock On stochastically complete manifolds.
\newblock {\em Soviet Math.Dokl}, 34, NO. 2:310--313, 1987.

\bibitem{Hsu89}
P.~Hsu.
\newblock Heat semigroup on a complete {R}iemannian manifold.
\newblock {\em Ann. of Prob.}, 17(3):1248--1254, 1989.

\bibitem{IC82II}
K.~Ichihara.
\newblock Curvature, geodesics and the {B}rownian motion on a {R}iemannian
  manifold {II}: Explosion properties.
\newblock {\em Nagoya Math. J.}, 87:115--125, 1982.

\bibitem{ITO50}
K.~It\^o.
\newblock On stochastic differential equations on a differentiable manifold.
\newblock {\em Nagoya Math. J.}, 1:35--47, 1950.

\bibitem{KA-LI}
L.~Karp and P.~Li.
\newblock The heat equation on complete {R}iemannian manifolds.
\newblock preprint.

\bibitem{LI89}
X.-M. Li.
\newblock Behaviour at infinity of stochastic solutions and diffusion
  semi-groups.
\newblock {\em Msc. dissertation}, 1989.

\bibitem{SC89}
L.~Schwartz.
\newblock Le semi-groupe d'une diffusion en liaison avec les trajectories.
\newblock In {\em Lecture Notes in Mathematics, 1372}, pages 326--342.
  Springer-Verlag, 1989.

\bibitem{TA91}
M.~Takeda.
\newblock On the conservation of the {B}rownian motion on a {R}iemannian
  manifold.
\newblock {\em Bull. London Math. Soc.}, 23:86--88, 1991.

\bibitem{Yau78}
S.~T. Yau.
\newblock On the heat kernel of a complete {R}iemannian manifold.
\newblock {\em J.de Math. pures Appl.}, 57(2):191--201, 1978.

\end{thebibliography}

\noindent
 Address: \\
\noindent  
Department of Mathematics, University of Warwick, Coventry CV4, 7AL, U.K.\\
\noindent
e-mail: xl@maths.warwick.ac.uk

\end{document}